\documentclass[10pt]{amsart}
\usepackage{amsmath,amssymb,latexsym,hyperref,url,graphicx}

\graphicspath{{figures/}}

\newtheorem{mystery}{Mystery}

\newcommand{\wc}{\square}
\newcommand{\bc}{\blacksquare}
\newcommand{\s}{\mspace{-2.76mu}}
\newcommand{\underdot}[1]{\text{\d{$#1$}}}

\begin{document}

\title{Regularity versus complexity \\ in the binary representation of $3^n$}
\author{Eric S. Rowland}
\address{
	Mathematics Department \\
	Tulane University \\
	New Orleans, LA 70118, USA
}
\date{October 31, 2009}

\begin{abstract}
We use the grid consisting of bits of $3^n$ to motivate the definition of $2$-adic numbers.  Specifically, we exhibit diagonal stripes in the bits of $3^{2^n}$, which turn out to be the first in an infinite sequence of such structures.  Our observations are explained by a $2$-adic power series, providing some regularity among the disorder in the bits of powers of $3$.  Generally, the base-$p$ representation of $k^{p^n}$ has these features.
\end{abstract}

\maketitle
\markboth{Eric Rowland}{Regularity versus complexity in the binary representation of $3^n$}

\section{Several mysteries}\label{mysteries}

The binary representation of a number $m$ can be thought of as encoding the unique set of distinct powers of $2$ that sum to $m$.  For example,
\[
	81 = 1010001_2 = 1 \cdot 2^6 + 0 \cdot 2^5 + 1 \cdot 2^4 + 0 \cdot 2^3 + 0 \cdot 2^2 + 0 \cdot 2^1 + 1 \cdot 2^0 = 2^6 + 2^4 + 2^0.
\]
We will display binary representations graphically by rendering $0$ and $1$ respectively as $\wc$ and $\bc$.  For reasons that will become clear, the convention in this paper when displaying binary representations graphically is to reverse the order of the digits relative to the standard ordering, so that higher indices are to the right.  For example, we write $81 = \underdot{\bc}\s\wc\s\wc\s\wc\s\bc\s\wc\s\bc = 2^0 + 2^4 + 2^6$, where the dot identifies the $2^0$ position (somewhat like a decimal point).

\begin{figure}
	\includegraphics{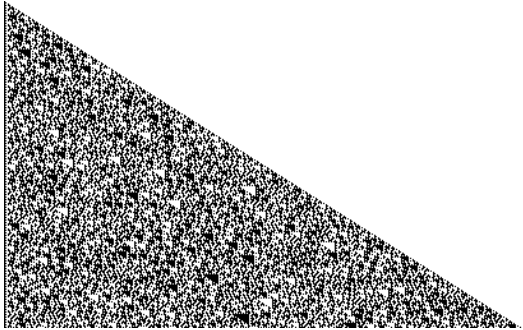}
	\caption{Powers of $3$ in base $2$.  The $n$th row consists of the binary digits of $3^n$, in order of increasing exponents.}
	\label{3n}
\end{figure}

The binary representations of the first several powers of $3$ grow steadily in length:
\begin{align*}
	3^0 &= \underdot{\bc} \\
	3^1 &= \underdot{\bc}\s\bc \\
	3^2 &= \underdot{\bc}\s\wc\s\wc\s\bc \\
	3^3 &= \underdot{\bc}\s\bc\s\wc\s\bc\s\bc \\
	3^4 &= \underdot{\bc}\s\wc\s\wc\s\wc\s\bc\s\wc\s\bc
\end{align*}
Figure~\ref{3n} displays a grid in which the $n$th row contains the binary digits of $3^n$.  Pictures like this were considered by Stephen Wolfram \cite[page 119]{NKS}.  Small triangles and other local features can be seen, but overall we get an impression of uniform disorder.  There is no global structure evident, aside from the right boundary of the pattern, which has slope $\log_2 3$.

There are some global regularities, however.  In particular, every column is eventually periodic.  This is because there are only finitely many (in fact $2^a$) states that can be assumed by the first $a$ columns taken together, so eventually the first $a$ columns return to a state that they have reached before, at which point they become periodic.

In fact, each column is not just eventually periodic but periodic from the start.  This is because each row has a unique predecessor, namely the integer obtained by dividing by $3$.  Put algebraically, $3$ is invertible modulo $2^a$ for every $a \geq 1$, so from a given row we may compute the previous row to as many bits as we want.

\begin{figure}[b]
	\includegraphics{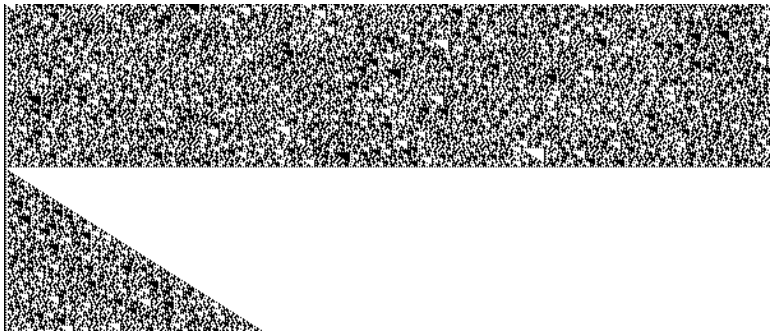}
	\caption{Part of the history obtained by periodically continuing each column up the page.}
	\label{3n-history}
\end{figure}

What if we try to compute $a$ bits of ``row $-1$'' --- the predecessor to the initial condition?  Certainly we can do this, and the result simply maintains the periodicity of the columns.  We can iteratively compute predecessors and thereby uniquely continue the picture up the page.  Figure~\ref{3n-history} shows the end of the unique infinite ``history'' leading up to the initial condition.

For those readers familiar with cellular automata, we mention that the ability to evolve the system backward in time is analogous to the same ability in a class of cellular automata whose local rules are bijective functions in the rightmost position.  As with multiplying each row by $3$ to form the next, information only propagates to the right and is not lost in these automata, and consequently they are reversible under the condition that the left half of each row is determined (say, all white) \cite{rowland06}.  Indeed, several themes of the present paper can be carried over to such cellular automata.

But what do the bits in such rows mean?  Rows $n < 0$ in the history illustrated in Figure~\ref{3n-history} do not represent integers, since they contain $1$s in positions arbitrarily far to the right:  The sum used to compute the value of such a row diverges.  For example, row $-1$ represents the ``infinite integer''
\begin{align*}
	\cdots\wc\s\wc\s\wc\s\wc\s\wc\s\underdot{\bc}\s\bc\s\wc\s\bc\s\wc\s\bc\s\wc\s\bc\s\wc\s\bc\s\wc\cdots
	&= 2^0 + 2^1 + 2^3 + 2^5 + 2^7 + 2^9 + \cdots \\
	&= 1 + \sum_{i=0}^\infty 2^{2i+1}.
\end{align*}
However, formally applying the geometric series formula to this divergent series produces
\[
	1 + 2 \sum_{i=0}^\infty 4^i \stackrel{?}{=} 1 + \frac{2}{1 - 4} = \frac{1}{3} = 3^{-1},
\]
which is certainly a natural object for row $n=-1$ to correspond to.  Similarly, row $-2$ represents
\begin{align*}
	\cdots\wc\s\wc\s\wc\s\wc\s\wc\s\underdot{\bc}\s\wc\s\wc\s\bc\s\bc\s\bc\s\wc\s\wc\s\wc\s\bc\s\bc\s\bc\s\wc\s\wc\s\wc\cdots
	&= 2^0 + \sum_{i=0}^\infty \left( 2^{6i+3} + 2^{6i+4} + 2^{6i+5} \right) \\
	&\stackrel{?}{=} 1 + \frac{2^3}{1 - 2^6} + \frac{2^4}{1 - 2^6} + \frac{2^5}{1 - 2^6} \\
	&= \frac{1}{9} = 3^{-2}.
\end{align*}
This is our first mystery, and in fact for each $n < 0$ there is a divergent series which produces $3^n$ under invalid applications of the geometric series formula.

\begin{mystery}\label{divergent}
Each power $3^n$ for $n < 0$ is the ``sum'' of a divergent series.
\end{mystery}

In order to resolve this mystery we must first encounter several additional mysteries --- all related to the first --- regarding the binary representation of $3^n$.

\begin{figure}
	\includegraphics{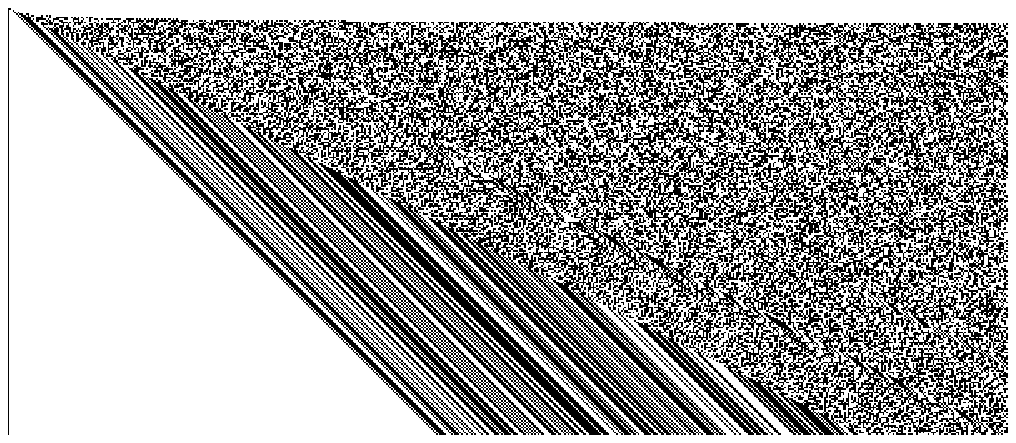}
	\caption{The subsequence $3^{2^n}$.  The $n$th row contains the first $600$ bits of $3^{2^n}$.}
	\label{32n}
\end{figure}

Since there are only two cell values ($\wc$ and $\bc$) in Figure~\ref{3n}, the period length of each column is a power of $2$.  (In fact for $a \geq 3$ the period length of $3^n \mod 2^a$ is $2^{a-2}$.)  Therefore, another consequence of the column periodicity is that row $2^n$ resembles the initial condition in several bits, since the periods of the first several columns will have just started over.  Brenton Bostick brought this ``local nestedness'' to my attention at the Midwest NKS Conference in 2005.  For example $3^{2^2} = 81 = \underdot{\bc}\s\wc\s\wc\s\wc\s\bc\s\wc\s\bc$ agrees with the initial condition $3^0 = \underdot{\bc}\s\wc\s\wc\s\wc\s\wc\s\wc\s\wc\cdots$ to three places to the right of the $2^0$ position.  Other terms in the subsequence $3^{2^n}$ agree with the initial condition to more places:
\begin{align*}
	3^{2^0} &= \underdot{\bc}\s\bc \\
	3^{2^1} &= \underdot{\bc}\s\wc\s\wc\s\bc \\
	3^{2^2} &= \underdot{\bc}\s\wc\s\wc\s\wc\s\bc\s\wc\s\bc \\
	3^{2^3} &= \underdot{\bc}\s\wc\s\wc\s\wc\s\wc\s\bc\s\wc\s\bc\s\bc\s\wc\s\wc\s\bc\s\bc \\
	3^{2^4} &= \underdot{\bc}\s\wc\s\wc\s\wc\s\wc\s\wc\s\bc\s\wc\s\bc\s\bc\s\bc\s\wc\s\bc\s\wc\s\bc\s\bc\s\wc\s\wc\s\wc\s\wc\s\bc\s\wc\s\wc\s\bc\s\wc\s\bc
\end{align*}
Figure~\ref{32n} shows many more rows, each row truncated at $600$ bits; the image can be efficiently produced with the following \emph{Mathematica} code.
\texttt{
\begin{tabbing}
	Arr\=ayPlot[ \\
		\> Rev\=erse /@ \\
			\> \> IntegerDigits[PowerMod[3, 2\^{}Range[0, 255], 2\^{}600], 2, 600] \\
	]
\end{tabbing}
}

The large triangular region of white cells indicates a sort of convergence to the initial condition:  The farther down we go in this image, the more columns have stabilized to the first bit in their period --- the bit in row $0$.  As Figure~\ref{32n} shows, each column (except the leftmost column) eventually becomes white, because row $0$ is simply $\cdots\wc\s\wc\s\wc\s\wc\s\wc\s\underdot{\bc}\s\wc\s\wc\s\wc\s\wc\s\wc\cdots = 1 = 3^0$.  In other words, $3^{2^n}$ ``converges'' bitwise to $1$ as $n \to \infty$.

This can be proven by Euler's theorem, which states that if $k$ is coprime to $b$ then $k^{\phi(b)} \equiv 1 \mod b$, where the Euler totient function $\phi(b)$ is the number of integers $1 \leq x \leq b$ that are relatively prime to $b$.  It is not difficult to convince oneself that if $p$ is prime then $\phi(p^{n+1}) = \frac{p-1}{p} \cdot p^{n+1} = (p-1) p^n$.  In our case, letting $k=3$ and $b=2^{n+1}$ gives $3^{2^n} \equiv 1 \mod 2^{n+1}$.  Letting $n \to \infty$ shows that every bit in $3^{2^n}$ eventually approaches the corresponding bit of $1$.

Perhaps we feel a little uneasy about giving much credence to this convergence, because certainly $3^{2^n}$ gets very large and far away from $1$ as $n$ gets large.  Thus we record it as another mystery.

\begin{mystery}\label{limit}
$\displaystyle{\lim_{n \to \infty} 3^{2^n} = 1}$.
\end{mystery}

\begin{figure}
	\includegraphics{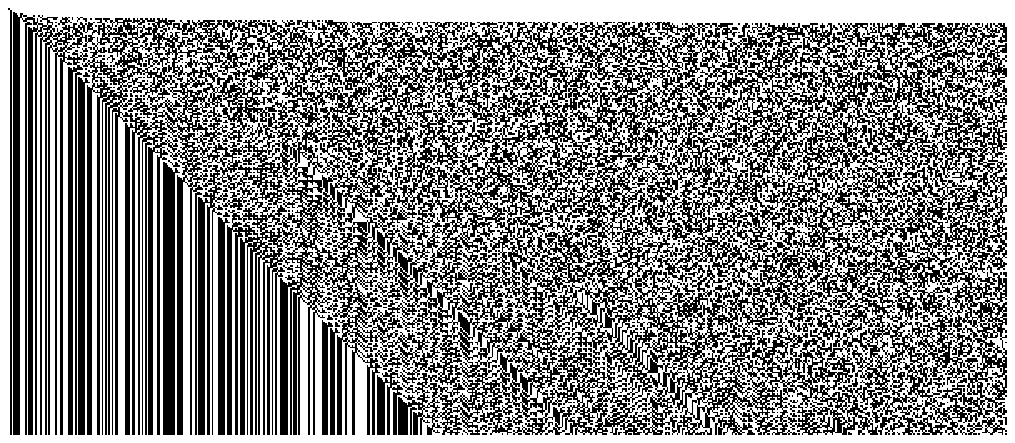}
	\caption{Bits of $3^{2^n}$, sheared so that the diagonal lines are now vertical.}
	\label{32n-sheared}
\end{figure}

In Figure~\ref{32n} we see additional structure as well --- surprising diagonal lines above the white triangular region.  More diagonals are filled in as we go down the page, so there appears to be another bitwise-convergent sequence here.  To change the diagonal lines into vertical lines, we make the first column white (for uniformity) and shear the image (shifting each row one position left relative to the row above it).  The result is shown as Figure~\ref{32n-sheared}.  Indeed these (shifted) rows are converging bitwise to something --- the row
\[
	c_1 = \cdots\wc\s\wc\s\wc\s\wc\s\wc\s\underdot{\wc}\s\wc\s\bc\s\wc\s\bc\s\bc\s\bc\s\bc\s\wc\s\wc\s\wc\s\bc\s\wc\s\bc\s\bc\s\bc\s\wc\s\bc\s\wc\s\wc\s\bc\s\wc\s\wc\s\bc\s\wc\s\bc\s\wc\s\wc\s\wc\s\bc\s\wc\s\wc\cdots.
\]
The shearing operation can also be effected by shifting the $n$th row left by $n$ bits; in other words, divide the $n$th row by $2^n$.  We may therefore write
\[
	c_1 = \lim_{n \to \infty} \frac{3^{2^n} - 1}{2^n}.
\]

\begin{figure}
	\includegraphics{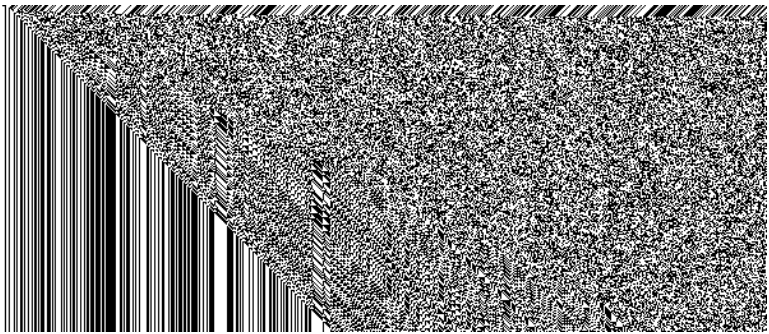}
	\caption{Bits of $\left((3^{2^n} - 1)/2^n - c_1\right)/2^n$, obtained from $3^{2^n}$ by twice subtracting the limit and shearing.}
	\label{32n-shearedagain}
\end{figure}

In Figure~\ref{32n-sheared} we also observe some secondary diagonal structures that were not easily visible before.  They are not as demarcated as the first set and seem to be interacting with the complex background.  In order to make the secondary diagonals vertical we would like to perform the same shearing operation.  However, first we need to subtract the limiting pattern $c_1$ from each row.  But subtract it how?  The limit is a divergent ``infinite integer'', but forming an integer from the first $a$ bits of $c_1$ and subtracting this integer from each row clears all the corresponding equal bits.  Once we have subtracted the limit, we divide by $2^n$ to remove the $n$ bits of $0$s on row $n$.  This produces Figure~\ref{32n-shearedagain}, in which the secondary diagonals are no longer muddied by the background but produce a clear limiting pattern themselves.  The new limit is
\begin{align*}
	c_2 &= \lim_{n \to \infty} \frac{\frac{3^{2^n} - 1}{2^n} - c_1}{2^n} \\
	&= \cdots\wc\s\wc\s\wc\s\wc\s\wc\s\underdot{\wc}\s\wc\s\wc\s\bc\s\wc\s\wc\s\bc\s\wc\s\wc\s\wc\s\bc\s\wc\s\bc\s\wc\s\wc\s\bc\s\bc\s\wc\s\bc\s\wc\s\wc\s\bc\s\wc\s\wc\s\bc\s\wc\s\bc\s\bc\s\wc\s\bc\s\wc\s\bc\cdots.
\end{align*}

It is natural to let $c_0 = \lim_{n \to \infty} 3^{2^n} = 1$ be the first limit.  If we continue to iterate this subtract-and-shear operation we continue to find convergent sequences of rows.  This means that, despite the apparent complexity in bits of $3^{2^n}$, every region can be decomposed into a sum of simple periodic regions.

The next limit
\[
	c_3 = \cdots\wc\s\wc\s\wc\s\wc\s\wc\s\underdot{\wc}\s\wc\s\wc\s\wc\s\wc\s\bc\s\bc\s\bc\s\wc\s\bc\s\bc\s\wc\s\wc\s\wc\s\wc\s\bc\s\bc\s\bc\s\wc\s\bc\s\wc\s\bc\s\bc\s\wc\s\wc\s\bc\s\wc\s\wc\s\bc\s\wc\s\wc\s\wc\cdots
\]
satisfies 
\[
	\lim_{n \to \infty} \frac{\frac{\frac{3^{2^n} - c_0}{2^n} - c_1}{2^n} - c_2}{2^n} - c_3 = 0.
\]
Let us take this expression and unravel it to see the structure better.  We find
\[
	3^{2^n} - (c_0 + c_1 2^n + c_2 2^{2n} + c_3 2^{3n}) \to 0
\]
as $n \to \infty$.  Replacing $2^n$ with $x$ reveals that this is a power series:
\[
	3^x = c_0 + c_1 x + c_2 x^2 + c_3 x^3 + O(x^4).
\]
Of course, we know a power series for $3^x$, namely
\[
	3^x = e^{x \log 3} = 1 + x \log 3 + \frac{1}{2!} (x \log 3)^2 + \frac{1}{3!} (x \log 3)^3 + \frac{1}{4!} (x \log 3)^4 + \cdots,
\]
so we might conjecture that $c_i = (\log 3)^i / i!$.

For $i=0$ we indeed have $(\log 3)^0 / 0! = 1 = c_0$.  But for $i=1$ the conjecture seems to fail, because $c_1$ is not a real number but an ``infinite integer''.  (In any case, the bits of $c_1$ don't resemble the binary representation of the real number $\log 3 = 1.00011001001111101010\cdots_2$.)

\begin{mystery}\label{log3}
``$\log 3$'' is not the real number $\log 3$.
\end{mystery}

\begin{figure}
	\includegraphics{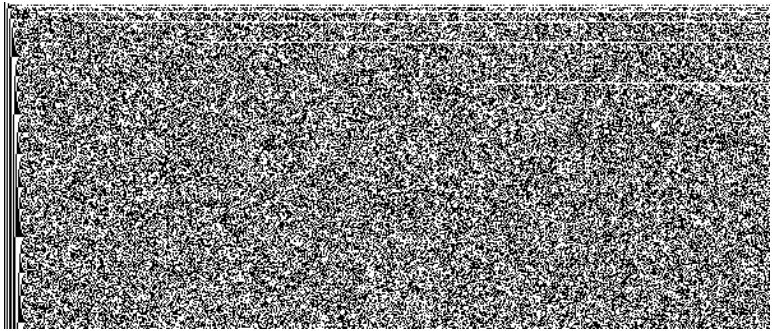}
	\caption{Binary representations of $(1 + 1/n)^n$ as real numbers, with most significant bits on the left.  The terms are slowly converging to $e = 10.10110111111000010101\cdots_2$.}
	\label{e}
\end{figure}

A final observation we make is that the direction of convergence is opposite that of real numbers.  All the sequences we have seen approach their limits by filling in bits from low indices to high indices, which is toward the right in our graphical convention of reversing the digits of integers.  A convergent sequence of real numbers, on the other hand, fills in bits from high indices to low indices.  Take the sequence $(1 + 1/n)^n$, for example.  Some terms of this sequence (as real numbers) are shown in Figure~\ref{e}.  The convergence proceeds from left to right, which is the same graphical direction but opposite numerical direction as the convergence of the sequence $3^{2^n}$ in Figure~\ref{32n}.

In the setting of bits of $3^{2^n}$, then, the low indices of a number are somehow stronger than the higher indices.  Therefore we should really think of the ``tail'' of numbers as being backward from the normal sense:  In this mode of convergence, two numbers are close to each other if their leftmost bits agree --- if their difference is divisible by a large power of $2$.  This is why we have chosen the convention that higher indices are to the right.

\begin{mystery}\label{close}
Two numbers are close if their difference is highly divisible by $2$.
\end{mystery}

\section{$2$-adic numbers}\label{2-adics}

Our four mysteries suggest that there is a notion of number presenting itself through the binary representation of $3^{2^n}$ that is quite different from the real numbers.  From Mystery~\ref{close} we must conclude that in some sense $2^i$ gets small as $i$ gets large, and the other mysteries support this conclusion.  Let us therefore make this a \emph{definition} instead of a mystery and introduce a new notion of ``size'' to make this precise.

Every rational number $r \neq 0$ has a representation $r = 2^\alpha \frac{n}{d}$ for integers $\alpha$, $n$, and $d$, where $n$ and $d$ are not divisible by $2$.  Moreover, $\alpha$ is unique.  We want $|r|_2$ to be large when $\alpha$ is small and small (but positive) when $\alpha$ is large.  A natural choice is to let $|r|_2 = 2^{-\alpha}$; this is called the \emph{$2$-adic norm} of $r$.  For example, $|64|_2 = 1/64$ and $|-691/2730|_2 = 2$.  Since $0$ is very highly divisible by $2$, let us define $|0|_2 = 0$.

Since large powers $2^i$ are small in the $2$-adic norm, a rational number can be a sum of arbitrarily large powers of $2$ when thought of $2$-adically, just as it can be a sum of arbitrarily large powers of $1/2$ when thought of as a real number.  For example,
\[
	\frac{1}{3} = 1 + \sum_{i=0}^\infty 2^{2i+1}.
\]

In fact, every rational number has the representation $\sum_{i=N}^\infty c_i 2^i$ for some integer $N$ and $c_i \in \{0,1\}$.  For example, $-1$ is rendered $2$-adically as
\[
	-1 = \frac{1}{1 - 2} = \sum_{i=0}^\infty 2^i = 2^0 + 2^1 + 2^2 + 2^3 + \cdots = \underdot{\bc}\s\bc\s\bc\s\bc\cdots,
\]
which illustrates that the $2$-adic representation of a number is really a ``limit'' of its representations modulo $2^a$ as $a \to \infty$.  For finite $a$, the representation modulo $2^a$ of course coincides with its two's complement representation.

In general, $a$ bits of the $2$-adic representation of a rational number $r = 2^\alpha \frac{n}{d}$ can be found by computing the inverse $d^{-1} \mod 2^a$, which is an integer, and multiplying by $2^\alpha n$.

One can check that the $2$-adic norm induces a metric $d(x,y) = |x-y|_2$ on the rational numbers, akin to the usual metric induced by the absolute value.  In particular, it satisfies the triangle inequality $|x-y|_2 + |y-z|_2 \geq |x-z|_2$.

There are some strange properties of this metric, however.  Perhaps the most immediate is that every triangle is isosceles:  If $|x-y|_2 = |y-z|_2$, then the triangle is isosceles by definition.  On the other hand, if $|x-y|_2 \neq |y-z|_2$, then
\[
	|x-z|_2 = \max(|x-y|_2, |y-z|_2).
\]
For example, if $x-y = 20$ and $y-z = 6$ then $|x-y|_2 = 1/4 \neq 1/2 = |y-z|_2$ and $|x-z|_2 = |26|_2 = 1/2$.

Naturally, we write $x_n \to x$ if the sequence of $2$-adic norms $|x - x_n|_2$ approaches $0$ as $n \to \infty$.  So indeed $3^{2^n} \to 1$ in the $2$-adic metric.

Of course, when we take a limit of rational numbers we may not get another rational number.  Traditionally, the real numbers can be constructed by taking limits of rationals with respect to the real metric; each real number has an expansion $\sum_{i=N}^\infty c_i 2^{-i}$ for $c_i \in \{0, 1\}$.  Similarly, we can take limits of rationals with respect to the $2$-adic metric and get a different completion of the rationals.  This completion is called the \emph{set of $2$-adic numbers}, and each $2$-adic number has a representation $\sum_{i=N}^\infty c_i 2^i$, where again $c_i \in \{0, 1\}$.  Like the real numbers, this set is complete --- it contains all its limit points.

It turns out that our power series $3^x = \sum_{i=0}^\infty (\log 3)^i x^i / i!$ is correct, but it must be interpreted not as a real power series but as a $2$-adic power series.  This means that ``$\log 3$'' is not the real number $\log 3$ but the $2$-adic number $\log 3$.  How do we compute it?  The function $\log (1-x)$ has a $2$-adic power series that coincides with its real power series:
\[
	\log (1 - x) = -\sum_{i=1}^\infty \frac{x^i}{i}.
\]
Of course, in the real metric this power series diverges at $x = -2$, so it cannot be used to compute the real $\log 3$.  But $2$-adically this series converges at $x = -2$ to the $2$-adic $c_1 = \log 3$.  Similarly, $c_2 = (\log 3)^2 / 2!$, $c_3 = (\log 3)^3 / 3!$, and so on.

To be precise, one must of course establish the standard objects of calculus over the $2$-adic numbers --- derivatives, power series, tests for convergence, etc.  We do not undertake this task here but refer the reader to texts on the subject.  The book of Gouv\^ea \cite{Gouvea} serves as a solid introduction, and Koblitz \cite{koblitz} provides a more advanced treatment.

\section{Generalizations}\label{generalizations}

The results of the previous section can be generalized in several directions, and we discover that the power series structure we have seen is quite common.

\begin{figure}
	\includegraphics{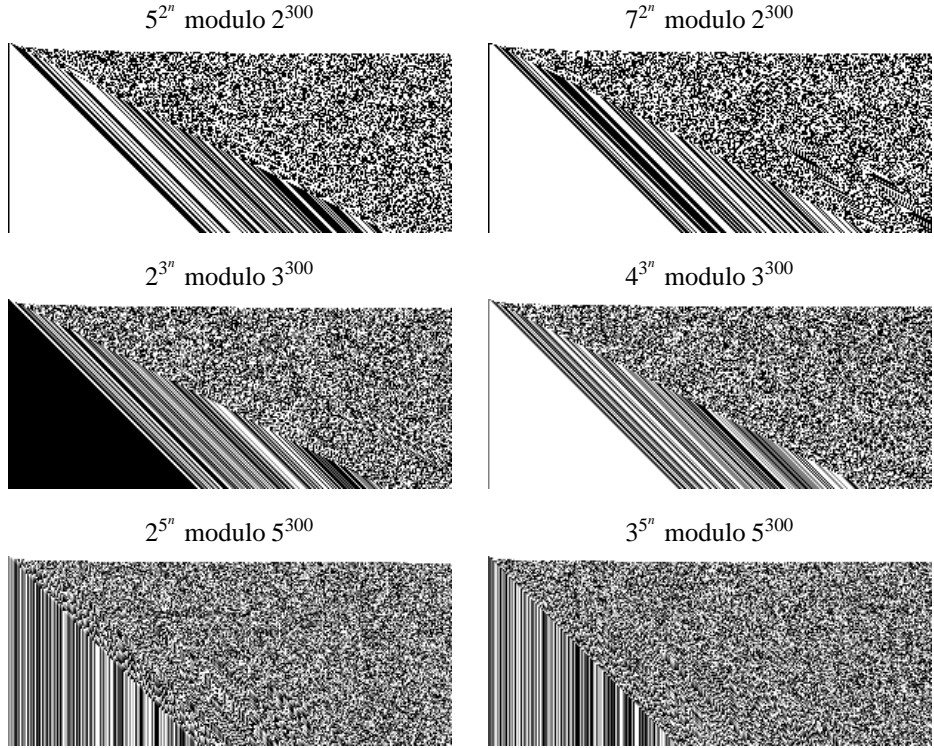}
	\caption{Powers $k^{p^n}$ in base $p$.}
	\label{otherpowers}
\end{figure}

Euler's theorem tells us that there is nothing particularly special about $3^{2^n}$, and in fact $5^{2^n}$ and $7^{2^n}$ have exactly analogous structures in binary, as shown in the first row of Figure~\ref{otherpowers}.  In general, if $k$ is odd then $k^x$ has a $2$-adic power series $1 + x \log k + \cdots$.  Gouv\^ea \cite[Section 4.5]{Gouvea} discusses the region of convergence of such power series.

What about other bases $b > 2$?  The second and third rows of Figure~\ref{otherpowers} show several examples.  To address these cases we briefly generalize the discussion to $p$-adic numbers for prime $p$.

Of course we may define $|x|_b$ for general $b$ (prime or composite) in the analogous way.  For primes $p$, $|x|_p$ is a norm on the set of rational numbers.  For composite $b$ it is not since in general $|x \cdot y|_b \neq |x|_b \cdot |y|_b$; for example, $|4|_4 = 1/4 \neq 1 = |2|_4 \cdot |2|_4$.  In fact, it is a theorem of Ostrowski that the $p$-adic norms and the real norm are (up to equivalence) the only nontrivial norms on the set of rational numbers.

Evidently $4^{3^n} \to 1$ in the $3$-adic metric.  However, the $3$-adic limit of $2^{3^n}$ is not $1$ but $2 + 2 \cdot 3^1 + 2 \cdot 3^2 + \cdots = -1$.  For general $k$ relatively prime to $p$, Euler's theorem provides that $k^{(p-1) p^n} \equiv 1 \mod p^{n+1}$.  In the limit, then, $k^{p^n}$ approaches a $(p-1)$th root of unity $1^{1/(p-1)}$.  This root of unity is congruent to $k$ modulo $p$ and is called the \emph{Teichm\"uller representative} of $k$.  This accounts for the vertical stripes in the base-$5$ digits of $2^{5^n}$ and $3^{5^n}$; the $5$-adic fourth roots of unity congruent to $2$ and $3$ modulo $5$ are irrational.  Note also that $2^{5^n} + 3^{5^n} \to 0$.

The $p$-adic power series of functions $f(x)$ other than $k^x$ are also evident in the base-$p$ digits of $f(p^n)$.  Let $F_n$, $C_n$, $M_n$, and $B_n$ be the sequences of Fibonacci, Catalan, Motzkin, and Bell numbers.  The sequences $C_{2^n}$ and $M_{2^n}$ have $2$-adic limits.  The sequences $F_{2^n}$ and $B_{2^n}$ do not have $2$-adic limits, but $F_{2^{2n}}$, $F_{2^{2n+1}}$, $B_{2^{2n}}$, and $B_{2^{2n+1}}$ do, giving some indication of the ubiquity of $p$-adic convergence in combinatorial sequences.

Finally, consider the factorial function $x!$.  The terms of the sequence $2^n!$, of course, become highly divisible by $2$, so $2^n! \to 0$ in the $2$-adic norm.  However, a theorem of Legendre implies that $|2^n!|_2 = 1/2^{2^n-1}$, and it turns out that
\[
	\frac{2^n!}{2^{2^n-1}}
\]
has a (nonzero) $2$-adic limit.

\end{document}